\documentclass[11pt]{article}
\usepackage{graphicx}
\usepackage[T1]{fontenc}
\usepackage[utf8]{inputenc}
\usepackage[dvips]{epsfig}
\usepackage{color}
\usepackage{amssymb}\usepackage{hyperref}\usepackage{amssymb}\usepackage{amscd}
\usepackage{float}
\usepackage{a4}

\newtheorem{theorem}{Theorem} \newtheorem{lemma}{Lemma}\newtheorem{remark}{Remark}\newtheorem{proposition}{Proposition}\newtheorem{Def}{Definition}
\newtheorem{corollary}{Corollary}\newtheorem{claim}{Claim}

\newcommand{\La}{\Lambda}
\newcommand{\R}{{\mathbb R}}   \newcommand{\N}{{\mathbb N}}
\newcommand{\T}{{\mathbb T}} \newcommand{\C}{{\mathbb C}}  \newcommand{\B}{\mathcal{B}}\newcommand{\PW}{\mathcal{PW}}
\newcommand{\bB}{{\mathbb B}}

\newcommand{\yy}{{\bf y}}\newcommand{\uu}{{\bf u}}\newcommand{\vv}{{\bf v}}\newcommand{\ww}{{\bf w}}\newcommand{\zz}{{\bf z}}
\newcommand{\xx}{{\bf x}}
\newcommand{\ee}{{\bf e}}\newcommand{\aaa}{{\bf a}}\newcommand{\sss}{{\bf s}}

\begin{document}
\author{Alexander~Rashkovskii \and Alexander~Ulanovskii  \and Ilya~Zlotnikov}
\title{On 2-dimensional mobile sampling}
\maketitle
\begin{abstract}
   Necessary and sufficient conditions are presented for several families of planar curves  to form  a set of stable sampling for the Bernstein space $\B_{\Omega}$ over a convex  set $\Omega \subset \mathbb{R}^2$. These conditions `essentially' describe  the mobile sampling property of these families for the    Paley-Wiener spaces $\PW^p_{\Omega},1\leq p<\infty$.
\end{abstract}

\section{Mobile Sampling Problem} The {\it classical sampling problem} is to determine  when every continuous signal (function) $f$ from a certain function space can be reconstructed from its  discrete samples $f(\lambda),\lambda\in\La$. The classical signal spaces are the Paley--Wiener spaces $\PW^p_\Omega$ of $L^p$-functions in $\R^d$ whose spectrum lies in a fixed set $\Omega\subset\R^d.$
When $p<\infty, $ the sampling problem asks for which discrete sets $\La\subset\R^d$ there exist positive constants $A,B$ such that
\begin{equation}\label{ss}
A\|f\|_p^p\leq\sum_{\lambda\in\La}|f(\lambda)|^p\leq B\|f\|_p^p,\quad \mbox{ for every } f\in \PW^p_\Omega.
\end{equation}

 A different method for the acquisition of  samples is when the samples of a multi-dimensional signal $f$ are taken by a mobile sensor that moves along a continuous path $\gamma$. The  {\it  mobile sampling problem } is then to reconstruct the signal from its samples on a continuous path or a union $P$ of continuous paths.
 In this case one needs to establish a `continuous variant' of the inequalities above:
 \begin{equation}\label{s}
A\|f\|_p^p\leq\int_P|f(u)|^p\,ds\leq B\|f\|_p^p,\quad \mbox{ for every } f\in \PW^p_\Omega,
\end{equation}where we assume that $P$ is locally rectifiable and integrate with respect to arc length.

 The mobile sampling problem has recently attracted much attention. We refer the reader to \cite{g,Jaming,u1,u2} for motivation and recent results.
 The sampling property in Paley--Wiener spaces of several families have been considered:

 (i) Parallel straight lines  in $\R^d$ (see i.e.   \cite{u1,u2,g}  and  references therein).

 (ii)  In  \cite{BW}, a sufficient condition for the Archimedes spiral is presented to form a set of stable sampling. In \cite{Jaming},  a wide family of  {\it spiraling curves} in $\R^2$ is introduced and  necessary and sufficient conditions for sampling in Paley--Wiener spaces with {\it convex symmetric spectrum} on these trajectories  obtained. 


In this paper, we consider the mobile sampling problem for three families of trajectories in $\R^2$. For each trajectory $P$ from one of these families, we present a necessary and sufficient condition for sampling in the Bernstein space $B_\Omega:=\PW_\Omega^\infty$ with {\it convex spectrum} $\Omega$. This condition `essentially' describes  the sampling property of $P$ for the Paley-Wiener spaces $\PW^p_\Omega$.

The rest of the paper is organized as follows. First, we give definitions
of the classical Paley-Wiener and Bernstein spaces, then a short list of notations.
In Section 4, we present the classical Beurling’s sampling theorem.
Our main results are formulated in Section 5 and proved in Sections 7-10.
In Section 6, we discuss the connection between the sampling in Bernstein
spaces and mobile sampling in Paley-Wiener spaces. In Section 11, we prove
some uniqueness theorems which may have independent interest. Finally, in
Section 12, we present some higher-dimensional results.

\section{Bernstein and Paley--Wiener spaces}

In what follows we will use the standard form of the Fourier transform:
\begin{equation}
    \hat f({\bf \yy}) = \int\limits_{\mathbb{R}^d} e^{-2 \pi i \langle \yy, \xx \rangle} f(\xx) d\xx,\quad \xx,\yy\in\R^d,
\end{equation}{}
where $\langle \cdot, \cdot \rangle$ is the usual inner product in $\R^d$.

We will consider the following classical  spaces of signals (functions):

\begin{Def}
Let $\Omega \subset \mathbb{R}^d,d\geq1,$ be a compact set.

{\rm 1}.  The Bernstein space $\B_{\Omega}$ consists of all continuous bounded functions in $\mathbb{R}^d$, which are the inverse Fourier transforms of tempered distributions supported by $\Omega$. Equipped with uniform norm $\|\cdot\|_\infty,$ $\B_{\Omega}$ is a Banach space.

{\rm 2}. Assume  $\Omega \subset \mathbb{R}^d$ has positive measure. The Paley-Wiener spaces $\PW^p_{\Omega},1\leq p<\infty,$ are defined as
$$
\PW^p_\Omega:=\B_\Omega\cap L^p(\R^d).
$$ Equipped with $L^p$-norm $\|\cdot\|_p,$ $\PW^p_{\Omega}$ is a Banach space.
\end{Def}

When $p=2$, the space  $\PW^2_\Omega $ is a Hilbert space consisting of all $L^2$-functions whose Fourier transform  vanishes a.e. outside $\Omega$.

Observe also that when  $\Omega\subset\R^d$ is a compact convex set, the space  $\B_\Omega$ admits an analytic description: It consists of all entire functions $f$ satisfying
$$
|f(\zz)|\leq C \exp\{-2\pi \max_{\uu\in\Omega}\langle \uu,\yy\rangle\}, \quad \zz=\xx+i\yy,\ \xx,\yy\in\R^d.
$$
\section{Notations}
Given $\vv\in \R^d,d\geq1,$ and $r>0,$ we denote by
$B_{r}(\vv):=\{\xx\in \R^d:|\xx-\vv|\leq r\}$   the closed  ball in $\R^d$ of radius $r$ centered at $\vv$. By $|E|$ we denote the ($p$-dimensional Lebesgue)  measure of a set $E\subset\R^d$ and  $\# E$ means the number of elements in $E$.

Set $\R^2_+:=\{\xx=(x_1,x_2):x_2\geq0\}$, $B^+_{r}(\xx):=B_r(\xx)\cap\R^2_+$ and $|\xx|:=\sqrt{x_1^2+x_2^2}$.

 Given sets $E,S\subset\R^2,$ $Q\subset\R$ and $\xx,\yy\in\R^2$, we write
$$
E+S=\{\xx+\yy: \xx\in E, \yy\in S\}, \ E-S=\{\xx-\yy: \xx\in E,\yy\in S\}, $$$$  QE=\{q\xx: \xx\in E, q\in Q\}, \ \
\mbox{dist}(\xx,S):=\inf_{\yy\in S}|\xx-\yy|.
$$
We say that the Hausdorff distance distance between $E$ and $S$ is $\leq\epsilon$ if $E\subset S+B_{\epsilon}(0)$ and $S\subset E+B_{\epsilon}(0)$.


\section{Sampling in $\B_\Omega$}
\begin{Def}
We say that a set $P\subset\R^d,d\geq1,$ is a sampling set (SS) for the Bernstein space $B_\Omega$, where $\Omega$ is a compact in $\R^d$, if there is a constant $C>0$ such that
\begin{equation}\label{sb}
\|f\|_\infty\leq C\|f|_P\|_\infty, \quad \mbox{for every } f\in \B_\Omega,
\end{equation}where$$
\|f|_P\|_\infty:=\sup_{\xx\in P}|f(\xx)|.$$
\end{Def}

\begin{Def}
{\rm 1}. A set $\La\subset\R^d,d\geq1,$ is called uniformly discrete (u.d.) if
\begin{equation}\label{delta}
\delta(\La):=\inf_{\lambda,\lambda'\in\La,\lambda\ne\lambda}|\lambda-\lambda'|>0.
\end{equation}The constant $\delta(\La)$ is called the separation constant for $\La$.

{\rm 2}. The lower uniform density of a set $\La\subset\R^d$ is defined as
$$
D^-(\La)=\lim_{r\to\infty}\inf_{\xx\in \R^d}\frac{\#\La\cap B_r(\xx)}{|B_r(\xx)|}.
$$
\end{Def}


In the classical situation where  $d=1$ and $\Omega$ is an interval  in $\R$, the sampling problem for $\B_\Omega$ was completely solved by Beurling:

\begin{theorem}\label{Beurling_one}{\rm (\cite{Beurling_one})}
Let $\Omega\subset\R$ be a compact interval. A set $P\subset\R$ is an SS for $\B_\Omega$ if and only if it contains a u.d. set $\La$ satisfying $D^-(\La)>|\Omega|.$
\end{theorem}

Observe that if $P\subset\R^d$ is an SS for $\B_\Omega$, then $P$ contains a discrete subset which is also an SS for $\B_\Omega$:

\begin{proposition}\label{r1}Assume $P\subset\R^d,d\geq1,$ is an SS for $\B_\Omega.$ Then there exists $\eta>0$ such that every  subset $\La\subset P$ satisfying
\begin{equation}\label{i}
P\subset \La+B_{\eta}(0)
\end{equation} is also an SS for $B_\Omega$.
\end{proposition}

We omit the proof, as it easily follows from   {\it Bernstein's inequality}, see \cite{OU}, p.21.

In particular, the set $\La$  in this result can be chosen to be u.d.


\section{Results}
In what follows, we assume that $\Omega\subset\R^2$ is a convex set of positive measure.

We will consider the  sampling problem for three families of curves in $\R^2$: parallel lines, dilations of a  convex closed curve around the origin and translations of a circle. For every family we present a sufficient and necessary condition for sampling in $\B_\Omega$.

\subsection{Parallel Lines}
Let $l\in \R^2$ be a straight line through the origin, and let $\vv_l$ be a unit vector orthogonal to $l$. Given any u.d. set $$H:=\{a_k\}_{k\in \mathbb{Z}}\subset\R,$$ consider the set of parallel lines
\begin{equation}\label{pl}P=l+        H\vv_l:=\bigcup_{k\in \mathbb{Z}} (l+a_k \vv_l).\end{equation}

\begin{theorem}\label{t1}
The set $P$ in {\rm (\ref{pl})} is an SS for $\B_\Omega$ if and only if $$D^-(H)\vv_l\not\in \Omega-\Omega.$$
\end{theorem}

\subsection{Dilations of a  convex curve}
Let  $D\subset \R^2$ be a closed convex set of finite positive measure such  that $0\in$ Int$(D)$. Denote by $\partial D$ the boundary of $D,$ by Ext$(D)$ the closed set of extreme points of $D$ and by $$D^o:=\{\xx\in\R^2:\langle \xx, \yy\rangle\leq 1, \yy\in D\}$$ the polar set of $D$.

\begin{figure}[h!]
\includegraphics[height=6cm]{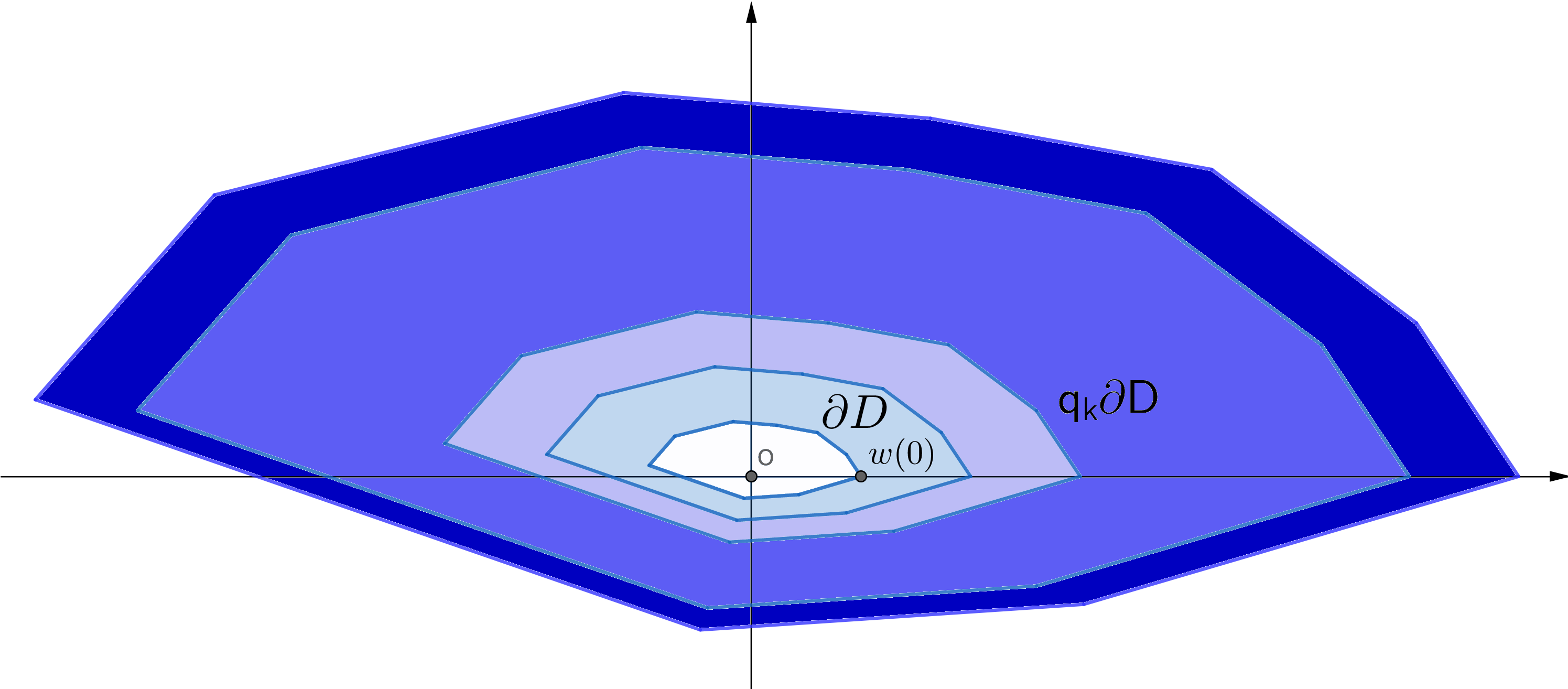}
\caption{Dilation of a convex curve}
\label{img1}
\end{figure}

Given a u.d. set $Q=\{q_k\}\subset (0,\infty)$, consider the set  \begin{equation}\label{dc}
P=Q\partial D:=\bigcup_{k=1}^\infty \bigcup_{\ww\in\partial D} \{q_k \ww\}.\end{equation}
Set $d^-(Q):=D^-(Q\cup(-Q))$.

\begin{theorem}\label{t2}
The set $P$ in {\rm (\ref{dc})} is an SS for $\B_\Omega$ if and only if \begin{equation}\label{the3}d^-(Q)\vv\not\in \Omega-\Omega,\quad \mbox{ for every }\vv\in \mbox{\rm Ext}(D^o).\end{equation}
\end{theorem}

The following is  simple corollary of this result and Remark \ref{r2} in sec. \ref{sec7}:

\begin{corollary}\label{co1}

{\rm (i)} Assume {\rm Ext}$(D)=D$. The set $P$ in {\rm (\ref{dc})} is an SS for $\B_\Omega$ if and only if $$  d^-(Q)\partial D^o\cap (\Omega-\Omega)=\emptyset.$$  In particular, if $D$ is the unit circle then  $P$ is an SS for $\B_\Omega$ if and only if {\rm Diam}$(\Omega)<d^-(Q)$.

{\rm (ii)} Let $D$ be the square $[-1,1]^2$. The set $P$ in (\ref{dc}) is an SS for $\B_\Omega$ if and only if $$(d^-(Q),0)\not\in\Omega-\Omega,\  (0,d^-(Q))\not\in\Omega-\Omega.$$
\end{corollary}

\subsection{Translations of  a Circle}
Take any circle  $T:=\{\xx\in\R^2:|\xx|=r\},r>0$. Let $V=\{\vv_k\}_{k=1}^\infty\subset\R^2$ be a u.d. set. Set
\begin{equation}\label{ts}
P=V+T:=\bigcup_{k=1}^\infty (\vv_k+T).
\end{equation}

\begin{theorem}\label{t3}
The set $P$ in (\ref{ts}) is an SS for $\B_\Omega$ if and only if $D^-(V)>0.$
\end{theorem}

\subsection{Remarks}
 Let  $\Omega\subset\R$ be a compact interval. Beurling's Theorem \ref{Beurling_one} solves the sampling problem for $\B_\Omega$ in terms of the lower uniform density $D^-(\La)$ of sampling set $\La$. The sampling property of $\La$ in the Paley-Wiener space $\PW^2_\Omega$ can be `essentially' described in terms of $D^-(\La)$, see i.e. \cite{OU}. See \cite{S} for necessary and sufficient conditions for sampling in $\PW^2_\Omega$. 

 If $\Omega\subset\R$ is  a disconnected set, already when it is a union of two intervals, the sampling property of u.d. sets $\La$ cannot be described   in terms of any density of $\La$. This is also the case for the spectra $\Omega\subset\R^d, d>1$, see \cite{OU, ou1}. However, the necessary density condition for sampling remains valid for general spectra $\Omega\subset\R^d$: Landau \cite{Landau} proved that if a u.d. set $\La$ is an SS for $\PW^2_\Omega$, then it satisfies  $D^-(\La)\geq |\Omega|$ (see \cite{no} for different simpler proof, which in particular  extends Landau's result to unbounded spectra).

 Observe that  for the mobile sampling (see definition in the next section)  there is no analogue of Landau's result.
 One can define a path density of trajectory $P$ as the `average length' covered by a curve.
     An example is constructed in \cite{g} showing that trajectories $P$ of arbitrarily small path density may nevertheless provide mobile sampling for $PW^p_\Omega$.
      Theorem \ref{tt} below reduces the mobile sampling problem for the Paley-Wiener spaces to the sampling problem for the Bernstein spaces. Hence, our Theorem \ref{t3}  presents another example in this direction. Observe also that  it easily follows from Corollary \ref{co1} (ii) and Theorem \ref{tt},  that a union of equidistant   squares $$P:=\bigcup_{n\in\N} \{\xx=(x_1,x_2):\max\{|x_1|,|x_2|\}=n\}$$  provides mobile sampling for $\PW^p_\Omega, 1\leq p<\infty$, for certain convex sets $\Omega$ of  arbitrarily large measure. However, by  Corollary \ref{co1} (i) and Theorem \ref{tt}, the union of equidistant circles
     $$
     P:=\bigcup_{n\in\N} \{\xx\in\R^2:|\xx|=n\}
     $$is not an SS for $\PW^p_\Omega$, whenever $\Omega$ is a convex compact set of measure $>\pi/4$.

\section{Sampling in $\B_\Omega$ versus Mobile Sampling in $PW^p_\Omega$}


Following \cite{g}, we call $P\subset\R^d$ a trajectory, if $P$ is a countable union of locally rectifiable (continuous) curves.

\begin{Def}
A trajectory $P\subset\R^d$ is called a stable sampling trajectory {\rm(}ST{\rm)} for $\PW^p_{\Omega},1\leq p<\infty,$ if   condition (\ref{s}) holds with some positive constants $A,B.$
\end{Def}{}

When the sampling set $\La\subset\R^d$ is discrete, it is well-known that the right inequality in  (\ref{ss}) is equivalent to the condition that there exist $r,C>0$ such that $\#(\La\cap B_r(\xx))<C$, for every $\xx\in\R^d$ (such sets $\La$ are called relatively uniformly discrete).

In the case of mobile sampling, one has
\begin{proposition}\label{p1}Let $1\leq p<\infty$,  $\Omega\subset\R^d$ be a compact set of positive measure and $P$ be a trajectory. The following conditions are equivalent:

(i) There is a constant $C$ such that
$$
\int_P |f(\uu)|^pds\leq C\|f\|_p^p,\quad \mbox{for every } f\in PW_\Omega^p;
$$

(ii) There are constants $r>0$ and $C>0$ such that
\begin{equation}\label{c}
 \sup_{\uu\in \R^d}\int_{P\cap B_r(\uu)}ds\leq C.
\end{equation}
\end{proposition}

We skip the proof which is similar to the corresponding proof for sampling on u.d. sets, see \cite{Ya}.



We will need a condition which prohibits   $P$ to  contain separated curves of arbitrarily small length:
\begin{equation}\label{cc}\mbox{{\sl There exists} }  r>0\mbox{  {\sl such that} }
 \inf_{\xx\in P}\int_{P\cap B_{\delta}(\xx)}ds\geq \delta,\quad \mbox{{\sl for every}  } \delta\leq r .
\end{equation}

The following result establishes a connection between sampling in Bernstein and mobile sampling in Paley--Wiener spaces:

\begin{theorem}\label{tt}
Let $1\leq p<\infty$ and $0<\epsilon<1.$  Let $\Omega\subset\R^d$ be a compact convex set of positive measure and  $P$ a trajectory satisfying (\ref{c}).

(i) If $P$ is an SS for $\B_\Omega$  and satisfies  (\ref{cc}), then it is an ST for  $\PW^p_{(1-\epsilon)\Omega}$.

(ii) If $P$ is not an SS for $\B_\Omega$, then it is not an ST for  $\PW^p_{(1+\epsilon)\Omega}$.
\end{theorem}

This is an analogue of the corresponding result for the discrete sampling  sets, see Theorem 5.30 in \cite{OU}. We omit the proof of Theorem \ref{tt} since it is rather similar to the proof of the mentioned result from \cite{OU}.

One can check that  part (i) of the theorem ceases to be true if condition  (\ref{cc})  is dropped.

\section{Proof of Theorem \ref{t1}}\label{sec7}

Before passing to the proof, we  recall several well-known facts.

 A sequence of u.d. sets $Q_k\subset\R,k\in\N,$ is said to converge weakly to a u.d. set $Q'$ if for every $R>0$ satisfying $\pm R\not\in Q'$, the Hausdorff distance between $Q_k\cap (-R,R)$ and $Q'\cap(-R,R)$ tends to zero as $k\to\infty$.

 The following lemma is well-known:
\begin{lemma}\label{l00}Assume u.d. sets $Q_k\subset\R$ satisfy $\inf_k\delta(Q_k)>0$, where $\delta(Q_k)$ is the separation constant defined in (5). Then there is a subsequence $Q_{k_n}$ which converges weakly to some  u.d. set $Q'$ satisfying  $D^-(Q')\geq \lim\sup_{n\to\infty}D^-(Q_{k_n})$.
\end{lemma}

One can also define the weak convergence of trajectories. However, in what follows we will  not  use this definition.

The next statement is obvious.
\begin{remark}\label{r2}
Assume $\Omega$ is a convex set and $\vv_{\theta}$ is a  non-zero vector with argument $\theta$, i.e. $\vv_{\theta} = |\vv_\theta|( \cos \theta, \sin \theta)$. Then $\Omega$ does not contain any segment of length $|\vv_\theta|$ parallel with $\vv_\theta$ if and only if $\vv_{\theta} \notin \Omega-\Omega$.
\end{remark}{}



Now, we pass to the proof of Theorem \ref{t1}. Without loss of generality we may assume that $l$ is a vertical line and $D^-(H)=1$.  To prove Theorem \ref{t1} we have to show that the set $P$ in (\ref{pl}) is an SS for $\B_\Omega$ if and only if \begin{equation}\label{line_main_condition} (1,0)\not\in \Omega-\Omega.\end{equation}


(i) Assume $P$ is not an SS for $B_{\Omega}$. We have to prove that (\ref{line_main_condition}) is not true.

Since $P$ is not an SS for $B_{\Omega}$, there is a sequence  of functions $f_k\in \B_{\Omega},k\in\N,$ such that
\begin{equation}
    \|f_k|_{P}\|_{\infty} \le \frac{1}{k}, \,\,\, \|f_k\|_{\infty} = 1 \mbox{   and   } \ f_k(\xx_k) \ge 1 - \frac{1}{k},
\end{equation}{}
for some points $\xx_k=(u_{k}, w_{k}) \in \mathbb{R}^2$.

Set
$$
g_k(\xx) := f_k(\xx+\xx_k).
$$
Clearly, we have $g_k(0) \ge 1 - \frac{1}{k}$ and $\|g_k|_{\{l + (H-u_k) (1,0)\}}\|_\infty \le \frac{1}{k}$.
By  Lemma~\ref{l00}, there is a subsequence $k_n$ such that the translates $H-u_k$ converge weakly to some u.d. set $H'$ satisfying $D^-(H')\geq D^-(H)=1$. By the compactness property of $\B_{\Omega}$ (see \cite{OU}), (taking if necessary a subsequence of $k_n$) we may assume that $g_{k_n}$ converge (uniformly on compacts in $\C^2$) to some non-trivial function $g \in \B_{\Omega}$. Clearly,
$\|g\|_\infty =g(0)= 1$  and  
\begin{equation}\label{ga} g(q,x_2)= 0,\quad q\in H', \ \ x_2\in\R.
\end{equation}{}

Fix any small $\varepsilon > 0$.
Since $D^-(H')\geq1$, by a result of Seip (see \cite{SeipRS}, Theorem 2.3), $H'$ contains a subset $H'' \subset H'$ such that the exponential system  $E(H''): = \{e^{2\pi i q t}:q\in H''\}$ forms a Riesz basis in $L^2(-(1-\varepsilon)/2 ,(1-\varepsilon)/2 )$.

Next, we invoke Pavlov's characterization of Riesz bases from \cite{Pavlov} to ensure that the  generating function of the exponential system above
\begin{equation}\label{genf}
\varphi(z_1) := \lim_{R\to\infty}\prod\limits_{q\in H'',|q|\leq R} \left(1 - \frac{z_1}{q}\right)
\end{equation}
is well defined, of exponential type $\pi (1 - \varepsilon)$, and  \begin{equation}\label{A2}|\varphi(x_1-i)|^2\in A_2,\end{equation} i.e. it belongs to the Muckenhoupt class $A_2$.

We briefly note that in  \cite{Pavlov} and \cite{SeipRS}, the zeroes of $\varphi$ are supposed to lie strictly above the real line.
However, one can overcome this obstacle by noting that the exponential system $E(H''+i):=\{e^{2\pi i (q+i)},q\in H''\}$
is also a Riesz basis for $L^2(-(1-\varepsilon)/2,(1-\varepsilon)/2)$, and its generating function $\tilde\varphi$ satisfies $\varphi(x_1-i)=\varphi(-i)\tilde\varphi(x_1)$. 

\begin{claim} There is a constant $\delta>0$ such that
\begin{equation}\label{Muck_est}
    |\varphi(x_1 - i)| \ge \frac{\delta}{1 + |x_1|^3},\quad x_1\in\R.
\end{equation} \end{claim}

This claim can be easily deduced from (\ref{A2}) and  Bernstein's inequality.

Recall that the generating function for a Riesz basis satisfies $|\varphi(x_1)|/(1+|x_1|)\in L^2(\R)$.  Choose any $q_0\in H''$ and set $\varphi_1(x_1):=\varphi(x_1)/(x_1-q_0).$

Then $\varphi_1\in \PW_{[-(1-\varepsilon)/2,(1-\varepsilon)/2]}$, and the points $\pm(1-\varepsilon)/2$ belong to the spectrum of $\varphi_1$.
Set
\begin{equation}
    \psi(z_1, z_2) := \frac{g(z_1, z_2) \sin^4\left(\frac{\varepsilon}{4} z_1\right)}{\varphi(z_1) {z_1}^4}.
\end{equation}
By (\ref{ga})  and (\ref{genf}),  $\psi$ is holomorphic in $\mathbb{C}^2$. By (\ref{Muck_est}), it belongs to $L^2$ on $(\R-i)\times\R$.
 Therefore, $\psi\in \PW_{\Omega'}$, where $\Omega'\subset\R^2$ is the spectrum of $\psi$.

 Consider the equality
\begin{equation}
    \varphi_1(z_1) \psi(z_1,z_2) = g(z_1,z_2) \psi_\varepsilon(z_1),\ \ \psi_\varepsilon(z_1):=\frac{\sin^4 \left(\frac{\varepsilon}{4} z_1 \right)}{z_1^4}.
\end{equation}
The spectrum of $\varphi_1$ contains the endpoints of the interval  $I:=[-(1-\varepsilon)/2,(1-\varepsilon)/2)]$ on the $x_1$-axis. The spectrum of $\psi_\varepsilon$  is the interval $I_\varepsilon:=[-\varepsilon,\varepsilon]$ on the $x_1$-axis.  One may now use an analogue of the Titchmarsh convolution theorem for higher dimensions:
$$\mbox{c.h.}(I+\Omega')=\mbox{c.h.(Sp}\,\varphi\cdot\psi)=\mbox{c.h.(Sp}\,g\cdot\psi_\varepsilon)\subset \mbox{c.h.}(I_\varepsilon+\Omega)\subset \Omega+B_{\varepsilon}(0),
$$where c.h. means the convex hull and Sp the spectrum. Clearly, the set  c.h.$(I+\Omega')$ contains a horizontal interval of length $|I|=1-\varepsilon$. This is also true for $\Omega+B_{\varepsilon}(0)$. It follows that $\Omega$ contains a horizontal interval of length $1-3\varepsilon.$ Using Remark \ref{r2}, we see that the point  $(1-3\varepsilon,0)\in \Omega-\Omega$. Since $\varepsilon$ can be chosen arbitrarily small, we conclude that (\ref{line_main_condition}) does not hold.

Note, that in this reasoning it is essential that $\Omega$ is a convex set.

(ii) Assume (\ref{line_main_condition}) does not hold. We have to show that $P$ is not an SS for $\B_\Omega$.
This is an easy consequence of  Beurling's Theorem 1.
Indeed, since translations of  $\Omega$ does not change the sampling property of $P$,  we may assume that $[-1/2, 1/2] \in \Omega$.
Using Theorem \ref{Beurling_one}, for every $\varepsilon>0$ there is a function  $f(x_1) \in \B_{[-1/2,1/2]}$ satisfying
$\|f|_{H}\|_\infty \le \varepsilon$ and $ \|f\|_{\infty} = 1.$ It follows that $P$ is not an SS for $\B_{[-1/2,1/2]}$. Therefore, $P$ is not an SS for $\B_\Omega$.

\section{ Auxiliary Results for the Proof of Theorem \ref{t2} }

 In this section,  $Q, D, \partial D$ have the same meaning as in Theorem \ref{t2}.

Denote by arg$\,\ww$ the argument of  vector $\ww$, i.e. the angle $ \theta$ such that $\ww=|\ww|(\cos\theta,\sin\theta) $.
We also denote by $\ww(\theta)\in\partial D$ the unique vector which lies on $\partial D$ satisfying arg$\,\ww(\theta)=\theta,$ $ -\pi<\theta\leq \pi$.

Recall that for every  convex or concave function $f(x)$ defined on an interval $I\subset\R$,  both one-sided derivatives of $f$  exist at every  interior point $w_0\in I$.
It follows that for every boundary point $\ww(\theta)\in \partial D$, both semi-tangent lines $\ww(\theta)+l_+(\theta)$  and $\ww(\theta)+l_-(\theta)$ exist, where $l_\pm(\theta)$ are straight lines through the origin.
In particular, if $\theta=0$ then there exist two lines $l_+(0)$ and $l_-(0)$ such that 
\begin{equation}\label{lp}
\mbox{dist}(\ww(\theta)-\ww(0), l_+(0)) = o(\theta),\quad  \theta\downarrow0,
\end{equation}
and $$
\mbox{dist}(\ww(\theta)-\ww(0), l_-(0)) = o(|\theta|),\quad  \theta\uparrow0.
$$

For the proof of Theorem \ref{t2}, we need two lemmas:

\begin{lemma}\label{l0}
Let $l_+(0)$ be the line satisfying (\ref{lp}). Assume a sequence of vectors $\xx_k$
satisfies  \begin{equation}\label{arg}\mbox{\rm arg}(\xx_k)>0,\quad k\in\N,\ \ |\xx_k|\to\infty, \ \ \mbox{\rm arg} (\xx_k)\to 0,\quad k\to\infty.\end{equation}
Then there exists a  subsequence $\xx_{k_n}$ and a u.d. set $Q'\subset\R$ such that \begin{equation}\label{udq} |\ww(\mbox{\rm arg}(\xx_{k_n}))|Q-|\xx_{k_n}|\, \mbox{converge weakly to}\, Q', \ D^-(Q')\geq \frac{d^-(Q)}{|\ww(0)|}.\end{equation}Condition (\ref{udq}) implies for every $R>0$ that
\begin{equation}\label{lp1}\mbox{\rm d}_R^+\left(Q\cdot \partial D-\xx_{k_n},\ww(0)+l_+(0)+Q'\cdot (1,0)\right)\to0,\ n\to\infty.\end{equation}  Here d$_R^+(A,B)$ denotes the Hausdorff distance between $A\cap B_R^+(0)$ and $B\cap B_R^+(0)$.
\end{lemma}

\begin{remark}

{\rm (i)} Note that if in   (\ref{arg}) we assume that the arguments of $\xx_k$ are negative and tend to zero, then  a similar to (\ref{lp1}) condition holds  with the `lower' semi-tangent line $l_-(0).$

{\rm (ii)}
Assume additionally that $l_+(0)=l_-(0)$, i.e.  there is a tangent  line to $D$ through $\ww(0)$. Then one may check that for every $R>0$ the Hausdorff distance between $$
(Q\cdot \partial D-\xx_{k_n})\cap B_R(0)$$  and  $$(\ww(0)+l_+(0)+Q'\cdot (1,0))\cap B_R(0)$$tends to zero as $n\to\infty.$

\end{remark}

\begin{lemma}\label{lt2}
Assume $P$ is not an SS for $\B_\Omega$. Then for every $n\in\N$ and $\epsilon>0$ there exist $f_n\in B_{(1+\epsilon)\Omega}$ and $\xx_n\in \R^2$
such that
\begin{equation}\label{nn}
|\xx_n|>n,\ \ \ \|f_n\|_\infty=1, \ \ \ \|f_n|_P\|_\infty<\frac{1}{n}, \ \ \ |f_n(\xx_n)|>1-\frac{1}{n}.
\end{equation}
\end{lemma}
\subsection{Proof of Lemma \ref{l0}}
Condition (\ref{udq}) follows easily from Lemma \ref{l00}. In what follows, for simplicity we assume that  \begin{equation}\label{trans}
 |\ww(\theta_n)|\cdot Q-|\xx_n|\mbox{  converge weakly to } Q',\quad  n\to\infty.\end{equation}

We have to deduce (\ref{lp1}) from (\ref{udq}). Before we proceed with the proof, observe that  (\ref{lp1}) is  intuitively clear. Indeed, since $\ww(0)+l_+(0)$ is a semi-tangent line to $\partial D$, when  $q_n$ and $|\xx_n|$ tend to infinity, the set  $(q_n\partial D-\xx_n)\cap B^+_R(0)$ is either empty or  `looks more and more like' a segment as $n\to\infty$.  However, the formal proof is somewhat technical.

Below we denote by $C$ different positive constants.

Recall that $D$ is a convex set of positive measure around the origin.

Given a straight line $l$ through the origin, denote by $\varphi(l),0\leq \varphi(l)<\pi,$ the  angle from the positive ray  $\R_+(1,0)$ to $l$ in the counterclockwise direction. Recall that $\ww(\theta)+l_+(\theta)$ and $\ww(\theta)+l_-(\theta)$ denote the semi-tangent lines to $\partial D$ at the boundary point $\ww(\theta)\in\partial D$, where we assume that $\varphi(l_+(\theta))\geq \varphi(l_-(\theta))$ for small positive values of $\theta$.

Clearly, for all small enough positive angles $\theta>\theta'$ we have
\begin{equation}\label{th}
\varphi(l_+(\theta'))\leq\mbox{\rm arg}\,(\ww(\theta')-\ww(\theta))\leq \varphi(l_-(\theta)),\quad 0<\theta'<\theta,
\end{equation}and
\begin{equation}\label{the}
|\ww(\theta')-\ww(\theta)|<C(\theta'-\theta),\quad 0<\theta'<\theta,
\end{equation}where one may take $C=2 |\ww(0)|$.

Clearly, $\varphi(l_\pm(\theta))\downarrow\varphi(l_+(0)),$
as $\theta\downarrow0$. Therefore, for every  $\epsilon>0$ there is an angle $\theta(\epsilon)>0$
such that \begin{equation}\label{th1}0<\varphi(l_+(\theta))-\varphi(l_+(0))<\epsilon, \quad 0<\theta<\theta(\epsilon).\end{equation}

 Set $\theta_n:=$arg$\,\xx_n$. By (\ref{arg}), $\theta_n>0$ and $\theta_n\to0$ as $n\to\infty.$ We assume that $n$ is large enough so that $\theta_n<\theta(\epsilon).$

Assume $\theta_n<\theta<\theta(\epsilon)$, and denote by $\ww_n(\theta)$ the point with argument $\theta$ lying on the semi-tangent line $\ww(\theta_n)+l_+(\theta_n)$:
\begin{equation}\label{wn}\ww_n(\theta)\in \ww(\theta_n)+l_+(\theta_n),\quad \mbox{\rm arg}\,\ww_n(\theta)=\theta. \end{equation}
From (\ref{th}), (\ref{the}) and (\ref{th1}) we may deduce that
\begin{equation}\label{then}
|\ww(\theta)-\ww_n(\theta)|<2C\epsilon(\theta-\theta_n),\quad \theta_n<\theta<\theta(\epsilon),
\end{equation}provided $\epsilon$ is sufficiently small.



Fix any $R>0$  satisfying  $\pm R\not\in Q'$. Then fix a positive number $\epsilon<1/R^2.$

To prove  (\ref{lp1}) we  show that  the Hausdorff distance between
\begin{equation}\label{set}
\left(Q\cdot\partial D-\xx_{n}\right)\cap B^+_R(0)
\end{equation}
and
$$
\left(\ww(0)+l_+(0)+Q'\cdot (0,1)\right)\cap B_R^+(0)
$$tends to zero as $n\to\infty.$
Since $\theta_n\to0$,   it suffices to check  this  for the Hausdorff distance between the set in (\ref{set}) and the set
$$
\left(\ww(\theta_n)+l_+(\theta_n)+Q'\cdot (\cos\theta_n,\sin\theta_n)\right)\cap B_R^+(0).
$$

Write $Q'\cap(-R,R)=\{q'(1),...,q'(m)\}$. By (\ref{trans}), for every large enough $n$,
$$
Q\cap\left(\frac{|\xx_n|-R}{|\ww(\theta_n)|},\frac{|\xx_n|+R}{|\ww(\theta_n)|}\right)=\{q_n(1),...,q_n(m)\},$$ where $$|\ww(\theta_n)| q_n(j)-|\xx_n|-q'(j)\to0,\quad n\to\infty,\quad j=1,...,m.
$$
Since arg$\,\ww(\theta_n)=$arg$\,\xx_n=\theta_n,$ this yields
\begin{equation}\label{theta}\left|\ww(\theta_n) q_n(j)-\xx_n-q'(j)(\cos\theta_n,\sin\theta_n)\right|\to0,\ n\to\infty,\ j=1,...,m.
\end{equation}

We see that it suffices to check that for every $j=1,...,m$, the Hausdorff distance between
\begin{equation}\label{sse}
\{q_n(j)\ww(\theta):\theta\geq\theta_n\}\cap (\xx_n+B^+_R(0))
\end{equation}and
\begin{equation}\label{sset}
(\xx_n+q'(j)(\cos\theta_n,\sin\theta_n)+l_+(\theta_n))\cap(\xx_n+ B^+_R(0))
\end{equation}tends to zero as $n\to\infty$.

Observe that for every sufficiently large $n$, condition $q\ww(\theta)\in \xx_n+B_R(0)$ implies $$\theta<\theta_n+2R/|\xx_n|, \quad q\in (C|\xx_n|/2, 2C|\xx_n|),$$ where we may take  $C=1/|\ww(0)|$. Hence, from (\ref{then}) one may easily check that the distance between the set $$\{q_n(j)\ww(\theta):\theta_n\leq\theta \leq \theta_n+2R/|\xx_n|\} $$ and the set\begin{equation}\label{sset1}\{q_n(j)\ww_n(\theta):\theta_n\leq\theta \leq \theta_n+2R/|\xx_n|\} \end{equation}
is less than $CR\epsilon<C\sqrt\epsilon$.

On the other hand, by (\ref{wn}), the point $q_n(j)\ww_n(\theta)$ has argument $\theta$ and lies on the line $q_n(j)\ww(\theta_n)+l_+(\theta_n).$
Let $\uu_j(\theta)$ be the point on $\xx_n+q'(j)(\cos\theta_n,\sin\theta_n)+l_+(\theta_n)$ satisfying arg$(\uu_j(\theta))=\theta$.
By (\ref{theta}), $|\uu_j(\theta)-q_n(j)\ww_n(\theta)|\to0$ as $n\to\infty$, which implies that for sufficiently large $n$,  the Hausdorff distance between the sets in (\ref{sse}) and (\ref{sset}) is less than $C\sqrt\epsilon.$ Since $\epsilon$ can be chosen arbitrarily small, this  proves (\ref{lp1}).

\subsection{Proof of Lemma \ref{lt2}} Since $P$ is not an SS for $\B_\Omega$, there is a sequence
of functions $g_k\in\B_\Omega$ satisfying
$$
\|g_k\|_\infty=1, \quad \|g_k|_P\|_\infty<\frac{1}{k}.
$$
For every $k$ choose a point $\yy_k$ such that $|g_k(\yy_k)|>1-1/k.$ If $$\lim\sup_{k\to\infty}|\yy_k|\to\infty,$$then condition (\ref{nn}) holds
for   $f_n(\xx):=g_{k_n}(\xx)\in \B_\Omega$, for some suitable subsequence $k_n$.

Assume  that the sequence $\yy_k$ is bounded. We may assume that it converges to some point $\yy_0\in \R^2.$ Using the compactness property of Bernstein spaces, see \cite{OU}, we may also assume that $g_n$ converge to some function $g_0\in \B_\Omega$. Clearly,  $g_0$ satisfies
$$
\|g_0\|_\infty=|g_0(\yy_0)|=1,\quad g_0|_P=0.
$$

Consider two cases.

1. Assume $g_0$ tends to zero fast in the sense that for every ${\bf m}\in(\N\cup\{0\})^2$ we have $|\xx^{\bf m}g_0(\xx)|\to0$ as $|\xx|\to\infty$, where $\xx^{\bf m}=x_1^{m_1}x_2^{m_2}, \xx=(x_1,x_2),{\bf m}=(m_1,m_2)$. Choose a point $\yy_{\bf m}$ such that
$$
\|\xx^{\bf m}g_0(\xx)\|_\infty=\max_{\xx\in\R^2}|\xx^{\bf m}g_0(\xx)|=|\yy_{\bf m}^{\bf m}g_0(\yy_{\bf m})|.
$$Clearly,
$$
|\yy_{\bf m}|\to\infty,\ \ \max_{\xx\in\R^2}|\xx^{\bf m}g_0(\xx)|\to\infty,\quad |{\bf m}|\to\infty.
$$Therefore, for a suitable subsequence ${\bf m}_n$,
the functions $$f_n(\xx):= \xx^{{\bf m}_n} g_0(\xx)/\|\xx^{{\bf m}_n} g_0(\xx)\|_\infty$$ belong to $\B_\Omega$ and satisfy condition (\ref{nn}).

2. If $g_0$ does not satisfy the decrease condition above, then there exist ${\bf m}\in(\N\cup\{0\})^2$ and a sequence of points $\yy_k$ such that
$$
|\yy_k|>k,\quad \quad |\yy_k^{\bf m} g_0(\yy_k)|\geq 1,\quad k\in\N.
$$

Consider the functions
$$
\varphi_k(\xx):=\xx^{\bf m} g_0(\xx)\mbox{sinc}^{2|{\bf m}|}(\epsilon(\xx-\yy_k)/2|{\bf m}|),
$$where
$$
\mbox{sinc}(\xx):=\frac{\sin x_1}{x_1}\frac{\sin x_2}{x_2},\quad \xx=(x_1,x_2).
$$Clearly, $\varphi_k$ belongs to $\B_{(1+\epsilon)\Omega}$,
vanishes on $P$ and $|\varphi_{k}(\yy_k)|\geq 1$. Moreover, for every $\xx$ in the disk $|\xx|<|\yy_k|/2$, we have
$$
|\varphi_{k}(\xx)|\leq \frac{|\xx^{{\bf m}}|}{\epsilon^{2|{\bf m}|}|\xx-\yy_k|^{2|{\bf m}|}}\leq
\frac{2^{|{\bf m}|}}{\epsilon^{2|{\bf m}|}|\yy_k|^{|{\bf m}|}}\to0,\quad k\to\infty.
$$Therefore, the functions
$$
f_n(\xx):=\frac{\varphi_{k_n}(\xx)}{\|\varphi_{k_n}\|}
$$satisfy (\ref{nn}) for a suitable subsequence $k_n$.

\section{Proof of Theorem \ref{t2}}

\subsection{Proof of Sufficiency} Assume  $P$ is not an SS for $\B_\Omega.$ We have to show that \begin{equation}\label{ssist}\mbox{There exists } \vv\in\mbox{\rm Ext}(D^o) \mbox{ satisfying } d^-(Q)\vv\in\Omega-\Omega.
\end{equation}

Fix any $\epsilon>0$. By Lemma \ref{lt2}, there exist  $f_n\in \B_{(1+\epsilon)\Omega}$  and $\xx_n\in \R^2$  satisfying (\ref{nn}).
Without loss of generality, we may assume that the sequence of arguments arg$(\xx_n)$ converges to zero. For convenience, we may assume that it converges to zero `from above', i.e. that it satisfies (\ref{arg}).

Using the compactness principle for Bernstein spaces, we may assume that
\begin{equation}\label{sf}  f_n(\xx+\xx_n)\mbox{ converge to some function } f\in \B_{(1+\epsilon)\Omega}.\end{equation}By convergence we mean the uniform convergence on compacts in $\C^2$. Clearly, the limit function $f$ satisfies $\|f\|_\infty=|f(0)|=1$.


 By  (\ref{lp1}), without loss of generality, we may assume that $f$ vanishes on the set of segments $(l_+(0)+\ww(0)+Q'(1,0))\cap B^+_R(0)$.
 Since $f$ is an entire function, it vanishes on the sets of lines $l_+(0)+\ww(0)+Q'(1,0)$.

 Denote by $\vv_l:=(\cos\varphi,\sin\varphi), |\varphi|<\pi/2,$ be a unite vector orthogonal to $l_+(0).$
 Denote by $A\subset\R$ the uniformly discrete set such that
 $$l_+(0)+\ww(0)+Q'(1,0)=l_+(0)+A\vv_l .$$It is easy to check that$$D^-(A)= \frac{D^-(Q')}{\cos\varphi}.$$
  Hence, by  (\ref{udq}), $$D^-(A)\geq \frac{d^-(Q)}{\cos\varphi |\ww(0)|}.$$
  Also, since  $f$ vanishes on $l_+(0)+A\vv_l$,  this set of lines is not an SS for $\B_{(1+\epsilon)\Omega}$. Then, by Theorem \ref{t1}, \begin{equation}\label{sist}\frac{d^-(Q)}{\cos\varphi|\ww(0)|}\vv_l\in (1+\epsilon)\Omega-(1+\epsilon)\Omega.\end{equation}

We will consider two cases:

Case 1. There is a unique point $\vv_0\in\partial D^o$ such that $\langle \ww(0),\vv_0\rangle=1.$

\begin{claim}$\vv_0\in\,${\rm Ext}$(D^o)$.\end{claim}

Indeed, assume $\vv_0\not\in\,$Ext$(D^o)$. Then $\vv_0$ is an inner point of a segment $I\subset D^o$. If this segment is  vertical, then  there are infinitely many points $\vv\in D^o$ satisfying $\langle \ww(0),\vv\rangle=1.$ If it is not vertical, we may find points $\vv\in I$ for which   $\langle \ww(0), \vv\rangle>1.$ None of the above is true. This shows that  $\vv_0\in\,$Ext$(D^o)$.

\begin{claim}\label{c2}  $\vv_0$ is orthogonal to $l_+(0).$ \end{claim}

Indeed, if not, then  there clearly exist infinitely many different vectors $\vv$ such that $\langle \vv,\ww(0)\rangle=1$ and $\langle \vv,\ww\rangle\leq 1$ for all $\ww\in D.$ But then $\vv\in D^o$, which contradicts the assumption above.

It follows that $\vv_0=|\vv_0|\vv_l$. 
Since  $\langle\vv_0, \ww(0)\rangle=1$, we conclude that $$\vv_0=\frac{\vv_l}{\cos\varphi|\ww(0)|}.$$Now, (\ref{ssist})
follows from (\ref{sist}), since it holds for every $\epsilon>0$.

Case 2. Assume there exist two points $\vv_1,\vv_2\in \partial D^o$ such that $\langle \ww(0),\vv_j\rangle=1, j=1,2.$ Then, clearly, the above holds for every
$\vv\in \partial D^o$ on the {\sl vertical} segment $I$ between $\vv_1$ and $\vv_2$. We may assume that $I$ is not a part of any larger segment which lies in $\partial D^o$.  Then, clearly,
$\vv_j\in$Ext$(D^o), j=1,2.$ Moreover, clearly, there is no point $\vv\in  \partial D^o\setminus I$ such that $\langle \ww(0) ,\vv\rangle=1.$

We may assume that $\vv_1$ lies `above' $\vv_2$. Similarly to Claim \ref{c2}, we show that $\vv_2$ is orthogonal to $l_+(0)$.
The rest of the proof repeats the proof above.

\subsection{Necessity}

For convenience, below we write $\B_\sigma:=\B_{[-\sigma,\sigma]}$.

We need a one-dimensional variant of Lemma \ref{lt2}.

\begin{lemma}\label{le11}Assume $\La\subset\R$ is a  uniformly discrete set, and let $\sigma:=D^-(\La)/2$.  Then for every $k\in \N$ there is a function $f_k\in \B_{\sigma+1/k}$ and a point $x(k)$ satisfying \begin{equation}\label{skkk}\|f_k\|_\infty=1, \ |f_k(x(k)|>1-1/k, \ |x(k)|>k, \ \|f_k|_\La\|_\infty<1/k.\end{equation}
\end{lemma}

Observe that by Theorem \ref{Beurling_one}, $\La$ is not an SS for $B_\sigma$.

The proof of Lemma \ref{le11} is similar to the proof of Lemma \ref{lt2}.

\begin{lemma}\label{le2} Assume a sequence of positive numbers $\sigma(k)$ converges to $\sigma>0$. Let $\La(k)\subset\R$ be u.d. sets satisfying $\inf_k \delta(\Lambda_k) > 0$ and such that   each $\La(k)$ is not an SS for $B_{\sigma(k)}$.  Then there is a sequence $x(k)$ with $|x(k)|\to\infty$ such that $\La(k)-x(k)$ converges weakly to some set $\La$ which is not an SS for $\B_\sigma$. If $\La(k)$ are symmetric, then $x(k)$ can be chosen positive.\end{lemma}

Proof. Use Lemma \ref{le11} to find  a sequence of functions
$f_k\in \B_{\sigma(k)+1/k}$ satisfying (\ref{skkk}) with $\La=\La(k)$. Clearly, if every $\La(k)$ is symmetric, we may assume $x_k>0$. Then the functions $f_k(x+x(k))$ converge to some non-trivial function $f\in \B_\sigma$.

By Lemma \ref{l00}, we may assume that the translates $\La_k-x(k)$ converge weakly to some set $\La$. Clearly, $f|_{\La}=0,$ which means that  $\La$ is not an SS for $\B_\sigma$.

\subsubsection{Proof of necessity}

 Assume $d^-(Q)\vv_0\in \Omega-\Omega,$ for some $\vv_0\in\,$Ext$(D^o)$.
We have to show that $P$ is not an SS for $\B_\Omega$, i.e. that for every small number $\eta>0$ there is a function $f_\eta\in \B_\Omega$ satisfying
\begin{equation}\label{fff}
\|f_\eta\|_\infty=1,\quad \|f_\eta|_P\|_\infty\leq\eta.
\end{equation}

Since $\vv_0\in\,$Ext$(D^o)$, there is a 
point $\ww_0\in D$, satisfying $\langle \vv_0,\ww_0\rangle=1.$ We may assume that $\ww_0=\ww(0)$, i.e. $\ww_0=(w,0), w>0$.
  Then the line  $\{\ww:\langle \ww,\vv_0\rangle=1\}$ is a semi-tangent line to $D$ at the point $\ww(0).$ We may assume that it is the `upper' semi-tangent line $\ww(0)+l_+(0).$ It will be convenient to write it in the form
  $$
  \ww(0)+l_+(0)=\ww(0)+ t\uu,\quad t\geq0,
  $$where $\uu$ is a unite vector parallel to $l_+(0)$ (and so, orthogonal to $\vv_0$).


Choose any positive sequence $ \theta_k\to 0,k\to\infty.$ Set
$$
\sigma(k):=\frac{d^-(Q)}{2|\ww(\theta_k)|}.
$$ Then $$\sigma(k)\to\sigma:=\frac{d^-(Q)}{2|\ww(0)|}.$$

By Theorem \ref{Beurling_one}, the symmetric set
$$
\frac{Q\cup(-Q)}{|\ww(\theta_n)|}
$$is not an SS for $\B_{\sigma(n)}$. Hence, by Lemma \ref{le2}, there exist $x(k)\to\infty $  such that the translates $Q|\ww(\theta_{k})|-x(k)$ converge weakly to some set $Q'$, which is not an SS for $\B_\sigma$.  By Lemma \ref{l0}, condition   (\ref{lp1}) holds with
$\xx_k:=x(k)(\cos\theta_k,\sin\theta_k)$.


Since $Q'$ is not an SS for $\B_\sigma$, it follows from the compactness principle for Bernstein spaces, that for every $\epsilon>0$ there exists $\delta>0$ such that there is a function $g(x)\in \B_{\sigma-\delta}$  satisfying
$$\|g\|_\infty=1,\quad \|g|_{Q'}\|_\infty\leq\epsilon.
$$

Consider the function $\varphi(\xx)$ defined as $$\varphi(\xx):=g\left(|\ww(0)|\langle \vv_0,\xx-\ww(0)\rangle\right).$$ It is easy to check that $\|\varphi\|_\infty=\|g\|_\infty=1$ and $$|\varphi(\xx)|\leq\epsilon,\quad \xx\in \ww(0)+t\uu+Q'(1,0),\ \ t\in\R.$$

By Remark \ref{r2}, condition $d^-(Q)\vv_0\in \Omega-\Omega$ means that $\Omega$ contains an interval $I$ parallel to $\vv_0$ and of length $d^-(Q)|\vv_0|.$ Since translations of $\Omega$ do not change the sampling property of $P$ for $\B_\Omega$, we may assume that $I $ is symmetric, $-I=I$. Observe that the spectrum of $g$ lies on $[-\sigma+\delta,\sigma-\delta]$, where $2\sigma=d^-(Q)/|\ww(0)|$. Then, clearly,  the spectrum of $\varphi $ lies on $(1-\delta)I$. Therefore, there is a small number $\delta'>0$ such that  $(1-\delta)I+B_{\delta'}(0)\subset\Omega.$

Now, choose a point $y_0$ such that $|g(y_0)|\geq 1/2.$ Then we have $|\varphi(\xx)|\geq 1/2$ for all $\xx$ on the line $\ww(0)+(y_0,0)+\R\uu.$


 Fix $R>0$ and consider the function
$$
\psi(\xx):=\varphi(\xx)\mbox{\rm sinc}(\delta'(\xx-\xx_0)),\quad  \xx_0:=\ww(0)-(y_0,0)-2R\uu.
$$
Then $\psi\in \B_\Omega$, $|\psi(\xx_0)|\geq 1/2$ and \begin{equation}\label{epsil}|\psi(\xx)|\leq\epsilon, \quad \xx\in \ww(0)+\R\uu+Q'(1,0).\end{equation}

Since $|$sinc$(\xx)|\to0$ as $|\xx|\to\infty,$  we may  assume that $R$ is so large that
\begin{equation}\label{epsi}|\psi(\xx)|\leq\epsilon, \quad |\xx-\xx_0|\geq R.\end{equation}

By Lemma 3, we may assume that the Hausdorff distance between the set $$(Q\partial D-\xx_k)\cap B^+_{4R}(0)$$ and $$(\ww(0)+l_+(0)+Q'(1,0))\cap B^+_{4R}(0)$$ tends to zero as $k\to\infty$.  We may also assume that $B_{R}(\xx_0)\subset B_{4R}^+(0)$. Then, the same is true for the sets $$(Q\partial D-\xx_k)\cap B_{R}(\xx_0)$$ and $$(\ww(0)+\R\uu+Q'(1,0))\cap B_{R}(\xx_0).$$
From (\ref{epsil}), by Bernstein's inequality, for all large enough $k$ we get
$$|\psi(\xx)|\leq C\epsilon,\quad \xx\in (Q\partial D-\xx_k)\cap B_{R}(\xx_0),$$where the constant $C$ depends only on the diameter of $\Omega$.

Finally, we see that the function
$$f(\xx):=\frac{\psi(\xx-\xx_k)}{\|\psi\|_\infty}$$ belongs to $\B_\Omega$ and satisfies
(\ref{fff}) with $\eta=C\epsilon$, where $\epsilon$ is any positive number.

\section{Proof of Theorem \ref{t3}}

 Recall that $P$ is defined in (\ref{ts}).

 (i) Assume that $D^-(V)=0$. We have to check that $P$ is not an SS for every space $\B_\Omega$, where $\Omega$ is a convex set of positive measure.

 The proof is easy. Indeed, from $D^-(V)=0$ it follows that there is a sequence of points $\xx_n$ such that the discs $B_{n}(\xx_n)$ do not intersect $V$. We may assume that $B_{\delta}(0)\subset\Omega$, for some $\delta>0$. Then the functions
  $$
  f_n(\xx):=\mbox{sinc}(\delta(\xx-\xx_n)), \quad n\in\N,
  $$belong to $\B_\Omega$ and satisfy $\|f_n\|_\infty=1$. It is obvious that  $$\|f_n|_P\|_\infty\leq \|f_n|_{\R^2\setminus B_{n}(\xx_n)}\|_\infty\to0,\quad n\to\infty,$$ which proves that $P$ is not an SS for $\B_\Omega$.

  (ii)
  Assume that $D^-(V)>0$. We have to check that $P$ is  an SS for every space $\B_\Omega$. The proof is a simple consequence of the uniqueness Theorem \ref{trci} below:
  Assume that $P$ is not an SS for $\B_\Omega$, i.e. there is a sequence of functions $f_n\in \B_\Omega$ satisfying
  $$
  \|f_n\|_\infty=1,\ \ \ \|f_n|_P\|_\infty\to0,\quad n\to\infty.
  $$

  Choose points $\xx_n$ such that $|f_n(\xx_n)|>1-1/n$ and set $g_n(\xx):=f_n(\xx+\xx_n).$ Then $$g_n\in \B_\Omega, \ \|g_n\|_\infty=1,\ |g_n(0)|>1/n.$$ Then a subsequence $g_{n_k}$ converges to some non-zero function $g\in \B_\Omega$.

  We may assume that the translates $V-\xx_{n_k}$ converge to some set $V'\subset\R^2$. We have $D^-(V')\geq D^-(V)>0$.
  It is clear that $g|_{V'+T}=0$. Theorem \ref{trci} yields $g=0$. Contradiction.


\section{Uniqueness sets}
Uniqueness sets play an important an important role in the sampling theory. In particular, Beurling \cite{Beurling_one}
proved that a u.d. set $\La$ is an SS for $\B_\sigma$ if and only if every weak limit of translates $\La-x_n$ is a uniqueness set for $\B_\sigma.$ A similar result holds in higher dimension. 

Below we consider subsets of $\R^2$ that are uniqueness sets for some classes of entire functions of exponential type and, in particular, of the Bernstein spaces. 
We believe such results are of independent interest.

Given en entire function $f$ in $\C^d$, let $Z_f=\{\zz\in\C^d:\: f(\zz)=0\}$ denote its zero set. For a generic $f$, the set $Z_f\cap\R^d$ is discrete and, as far as we know, only discrete (actually, u.d.) uniqueness sets $P\subset\R^d$ for entire functions of exponential type, as well as for $\B_\Omega$ with $\Omega=I_1\times\ldots \times I_d$  with $I_k=[-r_k, r_k]$, have been considered before (see, for example, \cite{Ber}, \cite{Ro}). Here we will be interested in the case of non-discrete uniqueness sets $P\subset\R^2$.

Note that, for any entire function $f$, the set $Z_f\cap\R^2$ is represented by the equations $f_\Re(\xx)=0$ and $f_\Im(\xx)=0$ with the real analytic functions $f_\Re={\rm Re}\, f$ and $f_\Im= {\rm Im}\, f$. Therefore, any set which is not a subset of a locally finite union of real analytic curves and discrete points in $\R^2$ is a uniqueness set for the whole class of entire functions in $\C^2$. In what follows, we work only with sets in $\R^2$ that are real analytic.

The idea of our considerations is, as in the case of discrete sets, getting control over the volume of the zero set $Z_f$ in the balls $\bB_t\subset\C^2$ in terms of a counting function, which in our case will be
\begin{equation}\label{theta1}\theta(t)=\# \{V\cap D_t\},\quad t>0,\end{equation}
where $V\subset\R^2$ is a discrete set related to $P$ and
 $D_t=\{\xx\in\R^2: \: |\xx|<t\}$ is an open disk in $\R^2$.

Three main ingredients are as follows. First, if a set $E\subset\R^2$ is non-discrete, then there exists at most one irreducible  analytic variety of complex dimension $1$ (i.e., analytic curve) in $\C^2$ containing $E$. Indeed, if $E\subset \chi_j$ for two irreducible analytic curves $\chi_1$ and $\chi_2$,  then $\dim_{\C} \chi_1\cap \chi_2=1$ and thus $\chi_1= \chi_2$.

Second, for any analytic curve $\chi$ and a point $\aaa\in\chi$, let $\sigma_{\chi,\aaa}(t)$ denote the volume of $\chi$ inside the ball
$$\bB_t(\aaa)=\{z\in\C^2:\: |\zz-\aaa|<t\}$$
and $\sigma_{\chi}(t)=\sigma_{\chi,0}(t)$. Then, by Lelong's bound for volumes of analytic sets, see \cite{LG}, Thm. 2.23,
\begin{equation}\label{LG}\sigma_{\chi,\aaa}(t)\ge \pi t^2\end{equation} for any $t>0$, with an equality if $\chi$ is a complex line.

Finally, we will use Jensen's formula for analytic functions in $\C^d$, see \cite{Ro_rev}:
$$ m_{\log|f|}(r):=\frac{1}{V_{d-1}}\int_{S_1}\log|f(r\zz)|\,dS_1(\zz)=\int_{r_0}^r \mu_f(t)\,t^{-2d+1}\,dt +C_{r_0,f},$$
where $dS_1$ is the normalized surface measure on $S_1=\partial\bB_1$, $V_{2d-2}$ is the volume of the unit ball in $\C^{d-1}$, and $\mu_f(t)$ is the volume, computed with the multiplicities, of $Z_f$ in $\bB_t$. This gives us for $d=2$
\begin{equation}\label{Jensen} m_{\log|f|}(r)\ge\frac1{\pi}\int_{r_0}^r \sigma_{Z_f}(t)\,t^{-3}\,dt + C_{r_0,f}.\end{equation}

Another, and more classical, form of Jensen's formula uses the intersections of $\chi$ with complex lines. Assume $0\not\in\chi$ and, for any point $\sss$ on the unit sphere $S_1$, let $n_{\sss\chi}(t)$ be the number of intersection points of $\chi$ with the line $\zz=\sss\zeta$, $\zeta\in\C$. Let
\begin{equation}\label{intersec} n_\chi(t)=\int_{S_1} n_{\sss\chi}(t)\,dS_1(\sss),\end{equation}
then
\begin{equation}\label{Jensen1} m_{\log|f|}(r)\ge\int_{r_0}^r \frac{n_{Z_f}(t)}t\,dt + C_{r_0,f}.\end{equation}

\subsection{Straight lines}

We start by considering the case when $P=\{l_k\}$ is a collection of straight lines in $\R^2$.
We will not assume, unlike in (\ref{pl}), that the lines are shifts of a single line $l$, the only condition being that none of the lines passes through the origin. Therefore, we can represent them as
\begin{equation}\label{flat}l_k=\{\xx\in \R^2:\: \langle \xx, \ee_k\rangle =1\}\end{equation}
for some vectors $\ee_k\in\R^2$. Denote
$\vv_k= \ee_k|\ee_k|^{-2}$ and let $\theta(t)$ be defined by (\ref{theta1}) for $V=\{\vv_k\}$. Note that if all $\ee_k=a_k^{-1}\vv_l$ for $a_k\in\R\setminus\{0\}$, this gives us precisely the set $P$ from (\ref{pl}).

\begin{theorem}\label{flatzeros} In the above setting, $P$ is the uniqueness set for entire functions of type $A$, provided
$$ \liminf_{t\to\infty}\frac{\theta(t)}{t}>\frac32\, A.$$
\end{theorem}

{\it Proof.} Let $L_k=\{\zz\in \C^2:\: \langle \zz, \ee_k\rangle =1\}$ be the complex lines containing $l_k$. Then any entire function $f\not\equiv 0$ vanishing on all $l_k$ also vanishes on all $L_k$, so
$$Z_f\supset Z:=\bigcup_k L_k.$$
By \cite{Papush},
\begin{equation}\label{papush}\sigma_{Z_f}(t)\ge \sigma_{Z}(t)=2\pi \int_0^ts\,n(s)\,ds,\end{equation}
where $n(s)$ is the amount of points $\vv_k$ inside the ball $\bB_s$. By the construction, all $\vv_k$  are real, so $n(s)=\theta(s)$ and
$$ \sigma_{Z}(t)=2\pi \int_0^ts\,\theta(s)\,ds.$$

There exist $A''>A'>A$  such that $\theta(s)\ge \frac32 A''s$ for all $s>r_0$ for some $r_0>0$. Therefore,
$$ \sigma_{Z}(t)\ge 2\pi \int_{r_0}^ts\,\theta(s)\,ds \ge  \pi A'' (t^3-r_0^3), \quad t>r_0,$$
and
$$
\frac1{\pi}\int_{r_0}^r \sigma_{Z_f}(t)\,t^{-3}\,dt  \ge  A''  \,(r-\frac32 r_0).
$$
If the function $f$ is of type $A>0$, then
$m_{\log|f|}(r)\le A'r$ for $r$ sufficiently big, which, by (\ref{Jensen}), is impossible.

\subsection{Dilations of circles}

Next, we will be concerned with dilations of the unit circle $\T=\{\xx\in\R^2:|\xx|=1\}$. Let $Q=\{q_k\}\subset \R_+$ be a discrete set, $\theta (t)$ be the counting function of $Q$, and
$$P=Q\,\T=\bigcup_k q_k \T.$$

\begin{theorem} $P$ is the uniqueness set for entire functions of type $A$ if
$$ \liminf_{t\to\infty}\frac{\theta(t)}{t}>\frac{A}{2\alpha},$$
where
$$\alpha=\int_{S_1} |s_1^2+s_2^2|^{1/2}\, dS_1(\sss).
$$
\end{theorem}

{\it Proof.}  Denote $\gamma=\{\zz\in\C^2:\: z_1^2+z_2^2=1\}$, the unique irreducible analytic curve in $\C^2$ containing the circle $\T$. For any $\sss\in S_1$ such that $s_1^2+s_2^2\neq 0$, the intersection $\sss q \,\gamma$ of the quadric $q\,\gamma$, $q>0$, with the line $\zz=\sss\,\C$ consists of two points given by $\zz=\sss\,\zeta$ with $(s_1^2+s_2^2)\zeta^2 =q^2$. Therefore,  (\ref{intersec}) gives us
$$n_{Q\gamma}(t)=2\int_{S_1} \theta(|s_1^2+s_2^2|^{1/2}t)\, dS_1(\sss).$$
Take any small $\epsilon>0$ and denote $E_\epsilon=\{\zz\in\C^2:\: |s_1^2+s_2^2|< \epsilon^2\}$. If $\theta(t)\ge \frac{A'}{2\alpha}t$ for some $A'>A\ge 0$ and all $t>\epsilon \,r_0$, then
$$n_{Q\gamma}(t)\ge 2\int_{S_1\setminus E_\epsilon} \theta(|s_1^2+s_2^2|^{1/2}t)\, dS_1(\sss) \ge (1-C_\epsilon)A' t$$
with $C_\epsilon\to 0$ as $\epsilon\to 0$,
and
$$\int_{r_0}^r \frac{n_{Q\gamma}(t)}t dt \ge (1-C_\epsilon)A'\,(r-r_0).$$
If an entire function $f\not\equiv0$ vanishes on $P$, then $Z_f\supset Q\gamma$. By (\ref{Jensen1}), it cannot have type
$A<A'$.

\subsection{Translations of circles}

Finally, we consider translations of a circle $T=\{\xx\in\R^2:|\xx|=r\}$.

\begin{theorem}\label{trci}
Let $\theta(t)$ be defined by (\ref{theta1}) for a discrete set $V\subset \R^2$. If
\begin{equation}\label{exptype}\lim_{t\to\infty}\frac{\theta(t)}{t}=\infty,\end{equation}
then $P=T+V$ is a uniqueness set for entire functions of exponential type. More precisely,
$P$ is a uniqueness set for entire functions of type $A>0$ if
\begin{equation}\label{exptypeA}B:=\liminf_{t\to\infty}\frac{\theta(t)}{t}> \frac3{2\sqrt2}A.\end{equation}
\end{theorem}

{\it Proof.} Denote $\gamma_0= \{\zz\in\C^2: \: \zz_1^2+\zz_2^2-r^2=0\}$, the unique  irreducible analytic curve containing $T$.

Let $f\not\equiv 0$ be en entire function in $\C^2$ of exponential type, vanishing on $P$. Then its zero set $Z_f$ contains $Z:=\cup_k \gamma_k$, where
$\gamma_k=\gamma_0+\vv _k$.
Since $\gamma_k\cap \gamma_j$ for any $j\neq k$ is either empty or a finite set, we have
$$ \sigma_{Z_f}(t)\ge \sigma_Z(t)=\sum_k \sigma_{\gamma_k}(t);$$
note that there is only finitely many $\gamma_k$ intersecting $\bB_t$.

 Take any $t>0$ sufficiently big and denote $K_t=\{k:\: \gamma_k\cap D_{t/2}\neq\emptyset\}$. Since $\bB_{t/2}(a)\subset \bB_t$ for any $a\in Z\cap D_{t/2}$, we get, by (\ref{LG}),
\begin{equation}\label{sigma_lbound} \sigma_Z(t)\ge\sum_{k\in K_t} \sigma_{\gamma_k}(t)\ge \theta (t/2)\,\pi (t/2)^2. \end{equation}

Assuming (\ref{exptype}), take any $N>0$ and let $r_0$ be such that
$\theta(t)> N t$ for all $t>r_0/2$.
 Then, by (\ref{sigma_lbound}), we get
$$ \frac1{\pi} \int_{r_0}^r \sigma_{Z_f}(t)\,t^{-3}\,dt \ge
\frac14\int_{r_0}^r \theta(t/2)\,t^{-1}\,dt
>  \frac{N}{8}\, (r -r_0),$$
which, in view of (\ref{Jensen}), means that $f$ cannot have finite type. This proves the first statement and a weaker version of the second one, with $B>8A$.

To prove the second statement in full, we need a sharper lower bound on the area of the variety $Z$ than (\ref{sigma_lbound}). Given $\epsilon\in (0,1)$,
let $r_0$ be such that
\begin{equation}\label{r_0}\theta(t)\,t^{-1}>(1-\epsilon)B\end{equation}
for all $t>r_0$, and take $\vv_k=(a_k,b_k)$ such that $r_0<|\vv_k|<t$. We have
\begin{equation}\label{sgk}
\sigma_{\gamma_k}(t)= \int_{\gamma_k\cap \bB_t}\,dm_2=t^2\int_{\gamma_{k,t}\cap \bB_1}\,dm_2=t^2\sigma_{\gamma_k,t}(1),
\end{equation}
where $\gamma_{k,t}$ is the analytic variety $\{(z_1-a_k/t)^2+(z_2-b_k/t)^2=t^{-2}\}$.

When $t\to\infty$, the varieties $\gamma_{k,t}$ converge to $\gamma_\infty:=\{z_1=\pm iz_2\}$. The convergence is not uniform in $k$, however we can choose $r_0$ such that, in addition to (\ref{r_0}), we have
\begin{equation}\label{gpm}\sigma_{\gamma_{k,t}}(1)\ge (1-\epsilon)\left(\sigma_{\Gamma_{k,t}^+}(1)+ \sigma_{\Gamma_{k,t}^-}(1)\right)
\end{equation}
for any $t\ge r_0$ and all $k$ with $|\vv_k|< t$, where $\Gamma_{k,t}^{\pm}$ are the complex lines
$$z_1-a_k/t=\pm i(z_2-b_k/t).$$
represent the families
$$\Gamma_t^{\pm}=\bigcup_k \Gamma_{k,t}^{\pm}$$
as in (\ref{flat}), by
$\langle (z_1,z_2), \ee_{k,t}^\pm\rangle =1$ with
$$\ee_{k,t}^\pm= \frac{(1,\pm i)}{c_k^\pm}\,t,$$
where $c_k^\pm=a_k\pm ib_k$.
The corresponding reference vectors $\vv_k$ are
$$\vv_{k,t}^\pm=\ee_{k,t}^\pm|\ee_{k,t}^\pm|^{-2}\frac{(1,\pm i)\overline{c_k^\pm}}{2t},$$
so $n_{\Gamma_{k,t}^{\pm}}(s)=\theta(\sqrt 2\,s\,t)$.

By (\ref{papush}), we have
$$\sigma_{\Gamma_{t}^{\pm}}(1)= 2\pi \int_0^1 s\, n_{\Gamma_{t}^{\pm}}(s)\,ds= 2\pi \int_0^1 s\, \theta(\sqrt 2\,s\,t)\,ds
= \pi t^{-2} \int_0^{\sqrt 2\, t} s\, \theta(\sqrt s)\,ds,$$
so  (\ref{sgk}), (\ref{gpm}) give us
$$
\sigma_Z(t)=t^2 \sum_k \sigma_{\gamma_{k,t}}(1)\ge (1-\epsilon)\,2\pi \int_0^{\sqrt 2\, t} s\, \theta(s)\,ds.
$$

Taking into account (\ref{r_0}), we get
\begin{eqnarray*} \frac1\pi\int_{r_0}^r \sigma_{Z_f}(t)\,t^{-3}\,dt &\ge& (1-\epsilon)\int_{r_0}^r\int_0^{\sqrt 2\, t} s\, \theta(s)\,ds
\,t^{-3}\,dt
\\
&\ge& (1-\epsilon)\int_{\sqrt2\,r_0}^{\sqrt 2\,r} \frac{\theta(s)}{2s} \left(2-\frac{s^2}{r^2}\right)\,ds\ge
(1-\epsilon)^2\frac{2\sqrt2}3\,B (r -3r_0),
\end{eqnarray*}
which contradicts  (\ref{Jensen}) if $f$ has type $A< \frac{2\sqrt2}3 \,B$. The proof is complete.

\bigskip
{\bf Remarks.} 1. As was mentioned before, functions $f$ from the Bernstein class $\B_\Omega$ have the bound
$$
|f(\xx+i\yy)|\leq C \exp\{2\pi H(\yy)\},
$$
where $H(\yy)=\max_{\uu\in\Omega}\langle \uu,\yy\rangle$ is the support function of $\Omega$, and so are of exponential type
$$A\le 2\pi\max_{|\yy|=1} H(\yy). $$
Actually, as follows from the Jensen's inequality, the constant $A$ in the both theorems can be chosen as
$$A=2\pi\int_{S_1} H(\yy)\,dS_1(\zz).$$

\medskip
2. The circle $T$ can be replaced with a trace of arbitrary irreducible entire analytic curve. Moreover, the whole collection $P$ can be formed by uniformly bounded arcs $T_k=\gamma_k\cap\R^2$ for irreducible analytic curves $\gamma_k$ satisfying $\dim\gamma_k\cap \gamma_j=0$, $j\neq k$, the set $V$ in the definition of the counting function $\theta(t)$ being formed by arbitrary points $\vv_k\in T_k$.

\section{Multi-dimensional Extensions}
Below we assume that the dimension $d\geq3$ and that $\Omega$ is a compact convex set of positive measure in $\R^d$.

1. The following extension of Theorem \ref{t1} holds true: {\it Let $l\subset\R^d $ be a hyperplane through the origin, $\uu_l$ be a unit vector orthogonal to $l$ and $H\subset\R$ be a u.d. set. Then $P=l+H\uu_l$ is an SS for $\B_\Omega$ if and only if $D^-(H)\uu_l\notin\Omega-\Omega$.} The proof is similar to the proof of Theorem \ref{t1}.

2. One may check that an analogue of Theorem \ref{t3} holds in higher dimensions for $P=V+T$, where $V\subset\R^d$ is a u.d. set and $T\subset\R^d$ is a $(d-1)$-dimensional sphere.

3. We guess that a multi-dimensional analogue of Theorem \ref{t2} is also true. In any case,  one may prove the following multi-dimensional variant of Corollary \ref{co1}:

(i) Let $Q\subset(0,\infty)$ be a u.d. set and $B_1(0)$ the unit ball in $\R^d$. The set $Q\partial B_1(0)$ is an SS for $\B_\Omega$ if and only if Diam$(\Omega)<d^-(Q).$

(ii) Let $D$ be a convex polytope around the origin. Then $Q\partial D$ is an SS for $\B_\Omega$ if and only if $d^-(Q)\vv\notin\Omega-\Omega$, for every vertex in the polar set (polytope) $\vv\in D^o.$

{}

\noindent A.R.: Stavanger University,  4036 Stavanger, Norway, alexander.rashkovskii@uis.no

\medskip

\noindent A.U.: Stavanger University,  4036 Stavanger, Norway, alexander.ulanovskii@uis.no

\medskip

\noindent I.Z.: Stavanger University,  4036 Stavanger, Norway, zlotnikk@rambler.ru  
\end{document}